\documentclass[11pt,reqno]{amsart}
\usepackage{graphicx}
\pdfoutput=1
\usepackage{amsfonts,amsmath,amssymb}
\usepackage{hyperref}
\usepackage{paralist}
\usepackage[linesnumbered,ruled,vlined,norelsize]{algorithm2e}

\def\1#1{\overline{#1}}
\def\2#1{\widetilde{#1}}
\def\3#1{\widehat{#1}}
\def\4#1{\mathbb{#1}}
\def\5#1{\frak{#1}}
\def\6#1{{\mathcal{#1}}}

\def\C{{\4C}}
\def\R{{\4R}}
\def\N{{\4N}}

\newcommand{\pz}{{P(z,\bar z)}}
\newcommand{\zz}{(z,\bar z)}
\newcommand{\La}{\Lambda}
\newcommand{\initial}{{\rm{in}\,}}

\newcommand{\initialterm}{{\rm{it}\,}}

\renewcommand{\Re}{\tmop{Re}}
\renewcommand{\Im}{\tmop{Im}}

\renewcommand {\a}{\alpha}
\renewcommand {\b}{\beta}

\usepackage{thmtools}
\declaretheoremstyle[bodyfont=\normalfont]{noncursive}
\declaretheorem{theorem}
\declaretheorem[numberwithin=section]{lemma}

\declaretheorem[style=noncursive,numberlike=lemma]{definition}

\declaretheorem[style=noncursive,numberlike=lemma]{remark}

\declaretheorem[style=noncursive,numberlike=lemma]{problem}
\renewcommand{\Re}{\mathop{\rm Re}\nolimits}
\renewcommand{\Im}{\mathop{\rm Im}\nolimits}
\newcommand{\im}{\ensuremath{\mbox{\rm Im}\,}}
\newcommand{\re}{\ensuremath{\mbox{\rm Re}\,}}

\newcommand{\CC}[1]{\mathbb{C}^{#1}}

\newcommand{\RR}[1]{\mathbb{R}^{#1}}
\newcommand{\dw}{\frac{\partial}{\partial w}}

\newcommand{\dz}{\frac{\partial}{\partial z}}

\newcommand{\la}{\lambda}

\numberwithin{equation}{section}

\newcommand{\aut}[1]{\mathfrak{aut}^{#1}}
\def\Label#1{\label{#1}}

\title[Normal forms in Cauchy-Riemann Geometry]{Normal forms in Cauchy-Riemann Geometry}

\author[M. Kolar]{Martin Kolar}
\address{Department of Mathematics and Statistics, Masaryk University, Brno}
\email{mkolar@math.muni.cz }
\author [I. Kossovskiy]{Ilya Kossovskiy*}

\address{\parbox{0.8\linewidth}{%
        Department of Mathematics, Federal University of Santa Catharina/\\ %
        Department of Mathematics and Statistics, Masaryk University, Brno}
    }
\email{ilyakos@gmail.com}
\thanks{*Supported by the Austrian Science Foundation.}
\author[D. Zaitsev]{Dmitri Zaitsev**}
\address{School of Mathematics, Trinity College, Dublin}
\email{zaitsev@maths.tcd.ie}
\thanks{**Supported in part by the Science Foundation Ireland grant 10/RFP/MTH2878.}

\begin{document}

\maketitle

\date{\today}

\begin{abstract}
One of effective ways to solve the equivalence problem and describe moduli spaces for real submanifolds in complex space is the normal form approach. In this survey, we outline some normal form constructions
in CR-geometry and formulate a number of open problems.
\end{abstract}

\tableofcontents

\section{Overview}

In the study of geometric structures on a manifold M, a normal form, corresponding to special
choices of coordinates adapted to the structure, is of fundamental importance. In case of
Riemannian metric, a possible choice of special coordinates is given by the normal coordinates given by
the exponential map. For a vector field $X$
that does not vanish at a point $p\in M$, special local
coordinates can be chosen in which this vector field is constant.
If, however, $X$ does vanish at $p$, its normal form, known as
{\em Poincare-Dulac normal form} \cite{ilyashenko}, exists in
general only in the formal sense, whereas its convergence is a
delicate issue depending on the presence of so-called small divisors.
Because of such a clear difference in behavior, a point where $X$
vanishes, is to be treated as a {\em singularity of $X$}.

The study of real submanifolds $M$ in complex spaces $\CC{n}$ is remarkable in that it exhibits both
regular and singular phenomena, with possible singularities being of very different nature. The
most obvious singularity for a real submanifold is a point
$p_0\in M$ for which the {\em complex tangent space}
$$T_p^{\CC{}}M:= T_pM \cap JT_pM$$
is of locally non-constant dimension for $p\in M$ near $p_0$. Such a point is called {\em CR singular}, and this property of the CR-structure at such a point is referred to as a {\em CR-singularity}. (Here $J$ is the
standard complex structure on $\CC{n}$.) For example, CR-singularities can occur for real surfaces in the complex space $\CC{2}$. Unless such a surface is embedded as a complex submanifold, it is locally biholomorphic at all CR nonsingular points to the unique flat normal form $\Pi=\{\im z_1=\im z_2=0\},\,(z_1,z_2)\in\CC{2}$. However, at a CR-singular point real surfaces in $\CC{2}$ possess biholomorphic invariants, as was shown in the celebrated paper of Moser and Webster \cite{mw} who constructed a normal form for  CR singularities of real-analytic surfaces in $\CC{2}$ satisfying an additional nondegeneracy assumption.  

If we now consider the class  {\em CR-submanifolds} (i.e., real submanifolds in $\CC{N}$ without CR-singular points), we discover that they still exhibit an extremely rich picture of possible singularities. The latter can viewed as points of ``non-uniformity'' of the CR-structure $(M,T^{\CC{}}M)$. One illustration
of possible singularities is as follows. For a germ of a real submanifold $M\subset\CC{N}$ at a point $p_0$ we consider local sections $X$ of the complex tangent bundle $T^{\CC{}}M\subset TM$ near $p_0$, and then for any $l\geq 1$  consider the subspace $g_l\subset  T^{\CC{}}_{p_0} M$   spanned by evaluating all possible  Lie brackets of order $\leq l$ of the above local sections at the point $p_0$. We then set 
$$\mathfrak g_1:= g_1,\quad \mathfrak g_l:=g_l/g_{l-1}.$$
The finite-dimensional graded Lie algebra  
\begin{equation}\Label{levitanaka}
\mathfrak g:=\oplus_{l\geq 1}\,\mathfrak g_l,
\end{equation}
where the Lie algebra operation is induced by taking Lie brackets,
is called the {\em  Levi-Tanaka algebra of $M$ at $p_0$}. The number and the dimensions of the graded components of  this graded algebra are immediate invariants associated with $(M,p_0)$,
as well as the Lie algebra structure.  
%If for some $l$ we have $\mathfrak g_l=T_{p_0}M$, then $p_0$ is called a {\em finite type point} (in the sense of Bloom-Graham \cite{bg}), otherwise $p_0$ is an {\em infinite type point}. 
Now a {\em singularity} of the CR-structure $(M,T^{\CC{}}M)$ can be viewed as a point, where the dimensions of the graded components of $\mathfrak g$ are not constant in any neighborhood. Further, even a point where these dimensions are locally constant but the natural map
\begin{equation}\Label{pgp}
p\mapsto \mathfrak g(p)\mapsto \mathcal M_s
\end{equation}
from $M$ into the moduli space of real Lie algebras of a fixed dimension $s$
is discontinuous, is still regarded as a CR-singular point for a given CR-structure. (The structure of the moduli space here should be understood here in the sense of the theory of invariants, see, e.g., \cite{vp}). A CR-structure without such singularities of its Levi-Tanaka algebra is called {\em regular} in what follows. Interestingly, not only the holomorphic classification of singularities of CR-structures but also the problem of holomorphic classification of {\em regular} CR-structures is  highly nontrivial and remains widely open  in its full generality. 

The key issue in relation to normal forms for real submanifolds is the problem of their {\em convergence}  or, more generally, the relation of the normal form to the holomorphic equivalence problem. Here one can distinguish between three possible scenarios. If a normal form exists merely in the formal sense (i.e\ all transformations and the normal form are given by formal power series, which are divergent in general, then such a normal form is good for identifying inequivalent  manifolds. It might still be possible that two formally equivalent CR-manifold are biholomorphically inequivalent.

 On the other hand, a convergent normal form gives a {\em complete} solution for the problem of equivalence of two objects.
%, and also allows for defining distinguished geometric objects on a manifold by transferring them back from a normal form where these objects have a particularly simple form. 
However, 
%CR-geometry exhibits a picture lying ``between'' the divergent and convergent  ones. That is, 
somewhat surprisingly, even a {\em divergent} formal normal form can in certain cases solve the {\em holomorphic} equivalence problem. This happens if real submanifolds of a certain class have the property that any {\em formal} map between {\em analytic} manifolds within the class is necessarily convergent. In this case, two real submanifolds are {\em holomorphically} equivalent if and only if some {\em formal} normal forms of them coincide. We refer to  Section 4 for examples of normal forms with this property, and also to the survey \cite{mirobzor} of Mir for the most up to date discussion of the relation between the formal and the holomorphic equivalences in CR-geometry. Some disadvantage of a normal form of  the latter kind  is that, unlike a convergent normal form, it does {\em not} allow to describe the moduli space of real-analytic CR-submanifolds under consideration, since it is not usually possible to identify among formal power series satisfying the normalization conditions the subset representing formal normal forms of real-analytic hypersurfaces.

The main goal of this paper is to give a survey of the normal form approach for the equivalence problem for CR-manifolds.
% (The CR-singular case is not addressed in this survey in detail).  
In what follows all CR-submanifolds in complex space $\CC{N}$ are assumed to be smooth and real-analytic,  and are all studied with respect to {\em local biholomorphic equivalences} $(M,p)\mapsto (M',p')$ of their germs at respective points $p,p'$. (By latter objects we mean local biholomorphic maps $F:\,(C^N,p)\mapsto (\CC{N},p')$ with $F(M)\subset M'$).

\section{Normal forms for Levi-nondegenerate CR-manifolds}

\subsection{Chern-Moser normal form}

The main example of a regular CR-structure is given by
{\em Levi-nondegenerate hypersurfaces}. The latter are real hypersurfaces for which  the {\em Levi form} 
\begin{equation}\Label{levi}
\6L_p\colon T_p^\C M \times T_p^\C M \to
\C\otimes (T_pM/ T_p^\C M), \quad \6L_p(X( p),Y( p)) = [X^{10},
Y^{01}]( p)
 \mod \C\otimes T_pM
\end{equation}
where $X$ and $Y$ are vector fields in $T^\C M$ and
$$X^{10}:=X-iJX, \quad X^{01}:= X +iJX,$$
are the corresponding $(1,0)$ and $(0,1)$ vector fields, is nondegenerate. In this case the map \eqref{pgp} is simply constant. A normal form for Levi-nondegenerate hypersurfaces was
constructed by Chern and Moser \cite{chern}. 
 A particular property of this normal form is its {\em convergence}. As a
consequence, geometry of the CR structure of a Levi-nondegenerate hypersurface can be studied
using its normal form. We outline below the result of Chern and Moser.  We now describe the master construction of Chern and Moser in detail, as it plays a fundamental role in the whole CR-geometry. The coordinates in $\CC{n+1},\,n\geq 1$ are denoted by $(z,w)=(z,u+iv)\in\CC{n}\times\CC{}$.

Probably, the crucial ingredient for having a formal normal form (the convergence is proved separately later) in Chern-Moser's case is to use the Poincar\'e local realization \cite{poincare} of a Levi-nondegenerate hypersurface:
\begin{equation}\Label{poincare}
v=\langle z,\bar z\rangle+O(3).
\end{equation}
Here $\langle z,\bar z\rangle$ is a non-degenerate Hermitian form on $\CC{n}$ which we may assume to have the form
\begin{equation}\Label{form}
\langle z,\bar z\rangle=\sum_{j=1}^n\varepsilon_j z_j\bar z_j,\quad \varepsilon_j=\pm 1,
\end{equation}
$O(3)$ is a power series in $z,\bar z,u$ each monomial of which has weight $\geq 3$, and the weights are assigned as:
\begin{equation}\Label{weights}
[z]=[\bar z]=1,\quad [u]=[v]=[w]=2.
\end{equation}
Poincar\'e realization immediately leads to the consideration of the {\em model} hypersurface $\mathcal Q$:
\begin{equation}\Label{model}
v=\langle z,\bar z\rangle,
\end{equation}
while all other Levi-nondegenerate hypersurfaces should be considered as perturbations of $\mathcal Q$. We call $\mathcal Q$ a{\em nondegenerate hyperquadric}, or shortly a quadric. An important object arising here is the (global) holomorphic automorphism group  of the quadric. The latter is a certain real form of the  projective group  $\mbox{Aut}\,(\mathbb{CP}^{n+1})$ and hence has the dimension $(n+2)^2-1$. The corresponding stability group of the origin $H:=\mbox{Aut}\,(\mathcal Q,0)$ has the dimension $n^2+2n+2$. 

It is crucial that Chern and Moser use the power series approach, in order to systematically analyse the {\em entire} freedom in perturbing a hypersurface near a distinguished point.   They suggest a procedure  to study  the action on germs of a hypersurface of the form \eqref{poincare} of the whole infinite-dimensional group of  {\em formal} biholomorphisms of $(\CC{n+1},0)$ preserving   \eqref{poincare}. They  set certain uniqueness conditions on a target hypersurface $M^*$ guaranteeing the uniqueness of a map sending a hypersurface $M$ of the form \eqref{poincare} into $M^*$. In view of the presence of the above group $\mathcal G$ preserving the quadric $\mathcal Q$, such a uniqueness of a map is possible only modulo the group $\mathcal G$ (note that the latter group preserves \eqref{poincare}). That is why for any (formal) map $F$ preserving \eqref{poincare} Chern and Moser first factor out an element $\psi\in \mathcal G$. That is, they present $F$ uniquely in the form 
$$F=G\circ\psi,\quad\psi\in \mathcal G,$$ where a map $G=(f,g),\,f(z,w)\in\CC{n},\,g(z,w)\in\CC{}$ has identity linear part and also satisfies 
$$\re g_{ww}=0.$$ 
We denote by $\mathcal F$ the space of maps normalized in the latter way. Now to find the desired uniqueness conditions on the target $M^*$ it is enough to study the action of the group $\mathcal F$ only. Then Chern and Moser analyse the transformation rule 
\begin{equation}\Label{tangency}
\im g(z,w)=\Phi^*(f(z,w),\overline{f(z,w)},\re g(z,w))|_{w=u+i\Phi(z,\bar z,u)}.
\end{equation}    
 Here $M=\{v=\Phi(z,\bar z,u)\},\,M^*=\{v=\Phi^*(z,\bar z,u)\},\,(f,g)\in\mathcal F:\,(M,0)\mapsto (M^*,0)$ (with $F=(f,g)$, again, being formal). Now one needs to use the weights \eqref{weights} to expand the map and the defining functions into sums of homogeneous polynomials: 
 \begin{equation}\Label{expand}
f=z+f_2+f_3+..., \quad g=w+g_3+g_4+...,\quad \Phi=\langle z,\bar z\rangle+\Phi_3+\Phi_4+...,
\end{equation}
and similarly for $\Phi^*$. Substituting \eqref{expand} into the transformation rule \eqref{tangency} and collecting for each fixed $m\geq 3$ all terms of weight $m$, Chern and Moser obtain the {\em homological equation}:
\begin{equation}\Label{homological}
\re(ig_{m}+\langle f_{m-1},\bar z\rangle)|_{w=u+i\langle z,\bar z\rangle}=\Phi_m-\Phi_m^*+...
\end{equation}
where dots stand for a finite number of terms depending on $\Phi_j,\Phi_j^*,f_{j-1},g_j$ (and their derivatives) with $j<m$ only. To analyse the possibility to solve the equations \eqref{homological} for the unknowns $(f_{m-1},g_m)\in\mathcal F$ on each fixed step $m$, we consider the linear operator
$$\mathcal L(f,g):=\re(ig+\langle f,\bar z\rangle)|_{w=u+i\langle z,\bar z\rangle}$$ mapping $\mathcal F$ into   the space $V$ of all (scalar) power series in $z,\bar z,u$ containing terms of weight $\geq 3$ only. $\mathcal L$ is called the {\em Chern-Moser operator}. We then consider some decomposition of the linear space $V$ of the form \begin{equation}\Label{decompose}
V=\mathcal L(\mathcal F)\oplus N.
\end{equation}
 Now it is nothing but a simple matter of linear algebra to see that {\em if the target defining function  $\Phi^*$ satisfies 
$$\Phi^*_m\in N\quad\forall\,m\geq 3,$$ then  solutions of the equations \eqref{homological} are parameterized by the kernel of the linear operator $\mathcal L$;  in particular, if the latter kernel is trivial, then  the equations \eqref{homological} determine the map $F=(f,g):\,(M,0)\mapsto (M^*,0)$ uniquely}. 

In this way, Chern and Moser reduce the hard problem of finding unique distinguished coordinates for a Levi-nondegenerate hypersurface to the, essentially, linear algebra problems of decomposing the space $V$ as $V=\mathcal L(\mathcal F)\oplus N$ and finding the kernel of $\mathcal L$. To describe the answer here, let us use the following notations: $\Phi_{kl}(z,\bar z,u)$ denotes a power series each term of which has the fixed total degree $k$ in $z$ and the fixed total degree $l$ in $\bar z$, and $\mbox{\bf tr}$ is the second order linear differential operator
$$\mbox{\bf tr}:=\sum\limits_{j=1}^{n}\varepsilon^j\frac{\partial^2}{\partial z_j\partial\overline{z_j}}$$
(here $\varepsilon_j$ are as in \eqref{form}). Now a smart power series computation due to Chern and Moser provides the following answer: 
$$\mbox{ker}\,\mathcal L=\{0\},$$
 and $N$ is the linear space of power series of the form 
\begin{equation}\Label{normalspace}
\sum_{k,l\geq 2}\Phi_{kl}(z,\bar z,u)\quad\mbox{with}\quad \mbox{\bf tr}\,\Phi_{22}=\mbox{\bf tr}^2\,\Phi_{23}=\mbox{\bf tr}^3\,\Phi_{33}=0.
\end{equation}
In fact, Chern and Moser prove the two latter statements at once by proving that {\em an equation $$\mathcal L(f,g)-\Psi\in\mathcal N,\,(f,g)\in\mathcal F,\,\Psi\in V$$ has a unique solution}. 

In fact, an interesting fact is obtained by computing the kernel of the operator $\mathcal L$ as defined on the space of all power series 
$$(f,g)=\left(\sum\nolimits_{k\geq 0}f_k,\sum\nolimits_{l\geq 0} g_l\right):$$ 
the kernel of such an operator is in one-to-one correspondence with the  Lie algebra $\mathfrak{hol}\,(\mathcal Q,0)$ of infinitesimal generators of the automorphism group  $\mbox{Aut}\,(\mathcal Q)$ (this algebra is called {\em the infinitesimal automorphism algebra of $(\mathcal Q,0)$}). Every element of this Lie algebra has the form 
$$X+\bar X,\quad X=f\dz+g\dw,$$
and $(f,g)$ appear to be precisely elements of the kernel. 
If we restrict to maps preserving the origin, then we obtain an isomorphism with the $5$-dimensional stability algebra 
\begin{equation}\Label{algebra}
\mathfrak h:=\mathfrak{aut}\,(\mathcal Q,0)\subset\mathfrak{hol}\,(\mathcal Q,0)
\end{equation}
 generating  the above group $H:=\mbox{Aut}\,(\mathcal Q,0)$. It is also not difficult then to verify from the equations \eqref{homological} and the Cramer rule that, in fact, each finite jet of a map $(f,g)$ bringing a hypersurface \eqref{poincare} to the normal form is rationally parameterized by elements of Lie algebra $\mathfrak h$. We summarize the results in the following

\begin{theorem}[Chern and Moser \cite{chern}]
For every real-analytic Levi-nondegenerate hypersurface $M\subset\CC{n+1},\,n\geq 1$ and any point $p$ in $M$ there exists a formal power series transformation bringing $(M,p)$ into normal form
\begin{equation}\Label{cmnf}
v=\langle z,\bar z\rangle+\Psi(z,\bar z,u),\quad \Psi(z,\bar z,u)\in N
\end{equation}
(where $N$ is the linear space of power series \eqref{normalspace}). The space of transformations bringing $(M,p)$ into  normal form is parameterized by the $5$-dimensional Lie algebra $\mathfrak h$, as in \eqref{algebra}. In fact, each finite jet of a normalizing transformation is parameterized by $\mathfrak h$ rationally. 
\end{theorem}

The remaining question now is the convergence of a normalizing transformation. The proof is based on 
the
presence of so-called {\em chains}, i.e.\ certain distinguished
real curves in the hypersurface $M$ that can be locally constructed as coordinate lines
\begin{equation}\Label{Gamma}
\{z=0,\quad v=0\}.
\end{equation} in
coordinates corresponding to the normal form. However, such a definition is based on the {\em a posteriori} established fact of convergence of a Chern-Moser normalizing transformation,  that is why Chern and Moser need a different definition which does not use the convergence property. Their original way to introduce chains is as follows. One chooses an arbitrary smooth real-analytic  curve $\gamma\subset M$, passing through the reference point $p$ and transverse to the complex tangent $T^{\CC{}}_p M$, and performs a transformation transforming $\gamma$ to \eqref{Gamma}. Then by a series of explicit (analytic) transformations  $(M,p)$ is brought to a form where all the normalization conditions in \eqref{normalspace} are satisfied except the trace conditions 
$\mbox{\bf tr}^2\,\Phi_{22}=\mbox{\bf tr}^2\,\Phi_{23}=\mbox{\bf tr}^2\,\Phi_{33}=0$. 
Inspecting the way $\mbox{\bf tr}^2\,\Phi_{23}$ depends on the initial curve $\gamma$, more precisely on its $2$-jet, Chern and Moser show that the requirement $\mbox{\bf tr}^2\,\Phi_{23}=0$ at each point amounts to a second order {\em analytic} ODE for the defining function of $\gamma$. The solution curves of this ODEs are called {\em chains} and are fundamental invariants of the underlying CR structure. 

The next step is to choose coordinates where a chain is a straight line. Then assuming the normalization conditions \eqref{normalspace} except the trace ones, it follows from the ODE that $\mbox{\bf tr}^2\,\Phi_{23}=0$ now holds identically along the chain. The remaining conditions $\mbox{\bf tr}^2\,\Phi_{22}=0$ and $\mbox{\bf tr}^2\,\Phi_{33}=0$ lead to decoupled systems of respectively $1$st and $3$rd oder ODEs along the chain, whose (simultaneous) solution finally yields the normal form. 

\begin{theorem}[Chern and Moser \cite{chern}]
A formal transformation bringing a real-analytic Levi-nondegenerate hypersurface $M\subset\CC{n+1},\,n\geq 1$ to the   normal form \eqref{cmnf} at a point $p\in M$ is convergent in some neighborhood of $p$. 
\end{theorem}
Theorem 2 implies, in particular, that {\em formal maps between two Levi-nondegenerate hypersurfaces are necessarily convergent}. 

In Section 3 below we explain a more geometric and less technical way to introduce chains, which somehow simplifies the quite involved original proof of Chern and Moser.

%Let us provide a small discussion of 
As Chern-Moser's normal form
%. This normal form is {\em complete}, to begin with. By that we mean that, among all normal forms in which  real hypersuadrics are assumed to be already in the normal form, the normal form of Chern and Moser has the property  that the number of parameters uniquely determining a normalizing transformation is smallest possible. Next, it 
is convergent, it yields the possibility to study the (local) complex geometry of  hypersurfaces. 
This normal form has the (convenient) property to be invariant under real shifts 
$$z\mapsto z,\quad w\mapsto w+a,\,a\in\RR{}$$
and complex scalings
$$z\mapsto\lambda z,\quad w\mapsto |\lambda|^2 w,\quad \lambda\in\CC{*}.$$ 
%Finally, it has the property that lowest possible degree terms in $z,\bar z$ are removed from the defining function in the normal form. 

 However, we shall emphasise that Chern-Moser's normal form is not anyhow a unique natural normal form. For some alternative normal form constructions for Levi-nondegenerate hypersurfaces we refer to the work \cite{zaitsevnf} of the third author.

\subsection{Beloshapka's theory of models}

Valeri Beloshapka in his research in 1980's - 2000's was developing a systematic way to extend Chern-Moser's strategy to as large as possible class of Levi-nondegenerate CR-submanifolds of high ($\geq 2$) codimension. This project was
% in a sense successfully 
completed in his 2000's paper \cite{universal}. The main goal of Beloshapka's research was finding CR-manifolds with maximal possible automorphism group. But as byproduct, Beloshapka developed a general strategy for constructing normal forms for Levi-nondegenerate CR-submanifolds.  We give below an overview of Beloshapka's work.

In what follows, we study a generic CR-submanifold $M\subset\CC{N}$ near the distinguished point $0\in M$. We split the ambient space as $\CC{N}=\CC{n}\times\CC{k}$, where $n:=\mbox{dim}\,T^{\CC{}}_0 M$ is the CR-dimension and $k$ is the real codimension of $M$ in $\CC{N}$.  The generic condition reads as $N=n+k$.   We denote the coordinates as $(z,w)\in\CC{n}\times\CC{k}$ and assume that the coordinates are chosen in such a way that $T^{\CC{}}_0 M=\{w=0\}$. We then call $M$ {\em Levi-nondegenerate}, if its Levi form defined as in \eqref{levi} (which is in our case valued in a $k$-dimensional real linear space) has zero kernel.

Let us start with the case $k\leq n^2$. 
\begin{definition}\Label{tot-nondeg} We call a  CR-submanifold $M$ as above with $k\leq n^2$ {\em totally nondegenerate}, if it is Levi-nondegenerate and its Levi-Tanaka algebra \eqref{levitanaka} has the form
$$\mathfrak g=\mathfrak g_1\oplus \mathfrak g_2$$
and, moreover, $\mbox{dim}\,\mathfrak g_2=k$. 
\end{definition}
We shall note that, in view of $\mbox{dim}\,M=2n+k$, {the integer $k$ is the maximal possible dimension of $\mathfrak g_2$}.  Thus the total nondegeneracy should be understood as the fact that {\em dimensions of the components of the Levi-Tanaka algebra of $M$ are maximal possible} (i.e., there is no linear dependence between the values at the origin of second order Lie brackets of vector fields spanning the complex tangents besides that given by purely algebraic reasons). Importantly, this condition is {\em generic}. 

For any totally nondegenerate CR-submanifold with $k\leq n^2$ Beloshapka proves \cite{belquadrics} that it can be locally represented in the Poincare form \eqref{poincare}, with the following correction: $\langle z,\bar z\rangle$ is now a $k$-valued hermitian form, such that its components  are linearly independent and the intersection of their kernels is zero. The weights are assigned as in \eqref{weights}. The corresponding model $\mathcal Q$, as in \eqref{model}, is also called {\em a nondegenerate quadric}. In relation to the initial germ $(M,0)$, it is called {\em the tangent quadric}. It is unique in any fixed coordinates.  Any {\em holomorphic} change of local coordinates amounts to a {\em linear} action on the defining function of the tangent quadric:
\begin{equation}\Label{transquadric}
\langle z,\bar z\rangle_2=\rho\cdot \langle Cz,\bar C \bar z\rangle_1,\quad \rho\in\mbox{GL}\,(k,\RR{}),\quad C\in\mbox{GL}\,(n,\CC{}).
\end{equation}
Here $\langle z,\bar z\rangle_1$ is the initial and $\langle z,\bar z\rangle_2$  the transformed vector-valued Hermitian form, respectively. Thus, one can consider the equivalence class of tangent quadrics with respect to the action \eqref{transquadric}. In fact, one can well define the moduli space $\mathcal M(n,k)$ of such quadrics \cite{moduli} and assign to the initial germ a unique point in $\mathcal M(n,k)$. Note that, for example, already the moduli space $\mathcal M(3,2)$ has positive dimension. Thus, a very important difference between the high-dimensional and the hypersurface cases is:

\smallskip

{\em In the case of higher codimension there may not be an unique tangent model; i.e., the equivalence class of tangent quadric depends in general on the base point $p\in M$.}

    \smallskip   

Note that quadrics (of high codimension) previously appeared in the work of Tanaka \cite{tanaka}, Naruki \cite{naruki}, and Kaup-Matsushima-Ochiai \cite{kmo}. They naturally occur as Shilov boundaries of Siegel domains.

Beloshapka shows \cite{belquadrics} that, similarly to the Chern-Moser case, for a (formal) map between two CR-manifolds of the form \cite{poincare} one can produce a homological equation similar to \eqref{homological}. The homological operator $\mathcal L$ is defined similarly as $$\mathcal L(f,g):=\re(ig+\langle f,\bar z\rangle)|_{w=u+i\langle z,\bar z\rangle},$$
however, in higher codimension case, it is vector-valued. Note that the existence of a formal map between CR-manifolds makes it possible to assume, after a linear coordinates change, that the tangent quadrics in the source and in the target are the same.  Now, 

\smallskip

{\em by finding a decomposition as in \eqref{decompose}, we obtain a complete (formal) normal form as in \eqref{cmnf} for the class of germs of submanifold with prescribed equivalence class of their tangent models; when fixing a tangent model $\mathcal Q$, the normalizing transformation is defined uniquely up to an element of the Lie algebra $\mathfrak h:=\mathfrak{aut}\,(\mathcal Q,0)$; moreover, each finite jet of the normalizing transformations depends on $\mathfrak h$ rationally}.

Remarkably, it is known that under Beloshapka's assumptions (which, in particular, imply the manifolds being $1$-nondegenerate and of finite type \cite{ber}) the formal equivalence of real submanifolds implies the holomorphic one (see \cite{mirobzor}), so that a formal normal form for those classes would always solve the holomorphic equivalence problem.

Another important issue here is the finite-dimensionality of the Lie algebra $\mathfrak h$. This fact was proved by Beloshapka using the general theory of linear differential operators. An alternative proof in more general cases was obtained by Baouendi, Ebenfelt and Rothschild \cite{ber98}.
%, whereas the structure of finite-dimensional Lie group on the germs of local automorphisms is due to the 3rd author \cite{z-germs}.

The above general strategy by Beloshapka leaves open the question on explicit description of the space $N$ in \eqref{decompose} for each fixed equivalence class of quadrics. In addition, it appears to be a quite subtle problem in general to describe the symmetries of a quadric, i.e. the algebra $\mathfrak h=\mathfrak{aut}\,(\mathcal Q,0)$ (note that the entire infinitesimal automorphism algebra can be always represented as $\mathfrak h_-\oplus\mathfrak h$, where $\mathfrak h_-$ is an affine algebra (i.e.\ consisting of affine transformations) acting transitively and freely on the quadric and given explicitly). A lot of work has been done on the latter question (see, e.g., the survey \cite{obzor} of Beloshapka on symmetries of real submanifolds). In particular, for $k=2$ and arbitrary $n$ a complete description is obtained by Shevchenko \cite{shevchenko}. Note that quadrics sometimes admit large automorphism groups and non-linear symmetries. 

The case of a general pair  $(n,k)$ was studied by Beloshapka in a subsequent series of publications (finalized by a construction of models for arbitrary $n,k>0$ in the paper \cite{universal}). For example, the next case to be studied after $k\leq n^2$ is the case $n^2<k\leq n^2(n+2)$. In this case, the total nondegeneracy amounts to the Levi-nondegeneracy of $M$ and to the conditions
$$\mathfrak g=\mathfrak g_1\oplus \mathfrak g_2\oplus \mathfrak g_3, \quad \mbox{dim}\,\mathfrak g_2=n^2,\quad  \mbox{dim}\,\mathfrak g_3=k-n^2$$
for the Levi-Tanaka algebra \eqref{levitanaka} (this conditions are, again, generic). The models corresponding to this case are {\em cubics}:
\begin{gather*}\left\{\begin{array}{lcl}
\mbox{Im}\,w_2=\langle z,\overline{z}\rangle\\
\mbox{Im}\,w_3=2\mbox{Re}\,\Phi(z,z,\overline z),\\
\end{array}
\right.\end{gather*}
where the coordinates in $\CC{N}$ are $(z,w_2,w_3)\in \CC{n}\times\CC{n^2}\times\CC{k}$, $\langle z,\overline{z}\rangle$ is a collection of $n^2$ linearly independent Hermitian forms (note that such a collection forms a basis so that we may assume it to be universal), and $\Phi(z,z,\overline z)$ is a collection of $k$ linearly independent cubic forms. Similarly to quadrics, cubics may have moduli space of positive dimension. Beloshapka shows that any totally nondegenerate manifold is a perturbation of a cubic:
\begin{gather*}\left\{\begin{array}{lcl}
\mbox{Im}\,w_2=\langle z,\overline{z}\rangle+O(3)\\
\mbox{Im}\,w_3=2\mbox{Re}\,\Phi(z,z,\overline z)+O(4),\\
\end{array}
\right.\end{gather*}
where the weights are chosen to be
$$[z]=1,\quad [w_2]=2, \quad [w_3]=3.$$
After that, a very similar strategy to the above one for quadrics of high codimension applies, including the above mentioned difficulties in determining the space $N$ explicitly and describing the stability algebra $\mathfrak h:=\mathfrak{aut}\,(\mathcal Q,0)$ of a cubic $\mathcal Q$. However, somewhat surprisingly, the cubic case appears to be much more rigid than the quadric case. Namely, as was shown by Gammel and the second author \cite{rigidity},

\smallskip

{\em the stability algebra of a cubic always consists of linear vector fields; furthermore, the stability group of a cubic is necessarily linear.}

\smallskip

As a corollary, one obtains:

\smallskip

{\em any map between two totally nondegenerate CR-submanifolds (see Definition~\ref{tot-nondeg}) is the case $n^2<k\leq n^2(n+2)$ is uniquely determined by the restriction of its differential onto the complex tangent $T^{\CC{}}_0 M$.
%in particular, this is the case for a transformation bringing such a CR-submanifold to Beloshapka's normal form.
}

\smallskip

The latter properties were extended by the second author for models of degree $4$ \cite{4degree}. 

After getting over certain technical difficulties, Beloshapka \cite{universal} was able to extend his theory of models to the case of arbitrary $n,k>0$. The definition of the total nondegeneracy works similarly to the quadric and cubic cases: we require that the CR-submanifold is Levi-nondegenerate and that the components of the Levi-Tanaka algebra have maximal possible dimensions (which is again a generic condition). We do not provide here the details of the structure of the algebra and that  of the models, since the latter are a bit technical. We only note that the length of the algebra, as well as the degree of the models, become arbitrarily large for fixed $n$ and large $k$.  We note also that it is still open whether there exist models of  degree $\geq 4$ with non-linear automorphisms.

The general strategy of constructing normal forms and the respective difficulties all remain the same in the general case. In certain cases Beloshapka's approach is realized explicitly, see, e.g., Ejov-Schmaltz \cite{es} and Beloshapka-Ejov-Schmaltz    \cite{bes2}. 

\section{Normal forms for Levi-degenerate hypersurfaces}

\subsection{The formal theory in $\CC{2}$} Levi degenerate hypersurfaces $M\subset\CC{2}$ were studied for the first time by J.J.Kohn.
in the context of the $\bar \partial$ - problem.
In \cite{kohn} he introduced the condition of finite type for hypersurfaces in $\CC{2}$,
as a main tool 
%lowest order invariant of a Levi degenerate hypersurface $M\subseteq \mathbb C^2$,  the
%type of the point. 
%The finite type condition 
%It was first introduced by J. J. Kohn in \cite{K},
 to study  subellipticity of the $\bar
\partial$ operator
on weakly pseudoconvex domains. For general $C^N$ the notion of finite type was extended later by Bloom-Graham \cite{bloom}.

Let us recall that in the original definition, the type of a point $p \in M$
 was introduced as the minimum number of commutators 
 of CR and anti-CR vector fields, needed to span the full complexified tangent space at the point (see \cite{kohn}).

Alternatively, the
type of $p$  can be defined as the maximum possible  order of contact
between $M$ and complex curves passing through $p$.

Let $r$ be a defining function for $M$ in a neighbourhood of $p$.
For a smooth
real valued function $f$ defined in neighbourhood of $0$ in
$\mathbb C$ let us denote by $\nu(f)$ the order of vanishing of $f$ at
$0$.

\begin{definition} $p$ is a point of finite type, if there
exists an integer $m$ such that
$$\nu (r\circ \gamma)\leq m$$ 
for all holomorphic maps $\gamma $  from a neighbourhood of $0\in
\mathbb C$ into $\mathbb C^2$, satisfying $\gamma(0)=p$ and
$\gamma^{\prime}(0)\neq 0$. The smallest such integer is called the
type of $p$.
\end{definition}

$M$ is Levi nondegenerate at $p$ if and only if $p$ is a point
of finite type two,  the case considered by the Chern-Moser theory. 

%Chern-Moser.

The first attempt to extend Chern-Moser's construction to Levi
degenerate manifolds is the work of Ph.~Wong \cite{wong}, who considered a special class of hypersurfaces of type
$4$ in $\mathbb C^2$, but did not succeed in constructing a complete normal form for this class. 
We also mention here the work of  Beloshapka \cite{belsmall} on the Levi-degenerate case,
and Ebenfelt \cite{ebenfeltjdg,ebenfeltC3} for $2$-nondegenerate hypersurfaces.
% Further results on the equivalence
%problem and normal form constructions were obtained
%in  \cite{S},
% (rigid hypersurfaces in $\mathbb C^2$),
%\cite{E},
%(a class of hypersurfaces of type three in $\mathbb C^{n}$, 
%\cite{BB}.
%(a class of hypersurfaces in $\mathbb C^2$ with a circular symmetry.

We now describe the complete formal normal form for finite type hypersurfaces, 
obtained by the first author in \cite{kolar}. 
It is not difficult to see that $p$ is a point of finite type $k$,
if and only if  there exist local holomorphic coordinates such that 
$M$ is
 defined
by
\begin{equation} v = P(z, \bar z) + o(\vert z \vert^k, u ),\ \label{2.9}\end{equation}
where the leading term  is a  nonzero real valued homogeneous
polynomial of degree $k$ without harmonic terms,
\begin{equation} P(z, \bar z) =  \sum_{j=1}^{k-1}
a_j z^j\bar z^{k-j}, \label{pzz}
\end{equation}
 with $a_j \in \mathbb C$ and $a_j =
\overline{a_{k-j}}.$ Note that \eqref{2.9} yields 
$$v=P(z,\bar z)+O(k+1),$$
if one assigns the weights
$$[z]=1,\quad [w]=k.$$
Thus, finite type hypersurfaces in $\CC{2}$ admit {\em models} of the kind 
$$\{v=P(z,\bar z)\},$$
where $P$ is as in \eqref{pzz}.

% The model hypersurface to $M$ at $p$ is
% 
% \begin{equation} M_H = \{(z,w) \in \cdva\ | \
%  v  = \sum_{j=1}^{k-1}
% a_j z^j\bar z^{k-j} \}. \label{2.10}\end{equation}

%we will write
%\begin{equation}S_k = \{(z,w) \in \mathbb C^2 \ \vert \   v = \ab k
%\}. \end{equation}
%Throughout this paper we consider the degenerate case, and assume
%that $k > 2$. 
The homogenous polynomial $P$ in \eqref{2.9} is not determined uniquely. 
We now introduce two
basic integer-valued invariants, needed for normalizing this polynomial leading term.
%Two basic integer valued invariants are used in the normal form
%construction in \cite{Ko1}. 
The first one is the essential
type of the model hypersurface, denoted  by $e$,
 which is the lowest index in (\ref{pzz}) for which $a_{e} \neq
0$.

For $e < \frac{k}2$,  the second invariant is defined as follows.
Let $e=m_0 <m_1 <\dots <m_s <\frac k2$ denote 
the indices in (\ref{pzz}) such that  $a_{m_i}\neq 0$. The invariant,
denoted by $d$, is the greatest common divisor of the numbers
$\ k-2m_0, k-2m_1,  \dots, k-2m_s$.

%The polynomial $P$ is not determined uniquely by 
%the form 
%(\ref{pzz}).
We need to impose two normalization
conditions which make the leading polynomial uniquely determined, namely
$$a_e = 1,$$ and

\begin{equation}
 \arg a_{m_{i+1}} \in [0,\frac{2\pi}{q_i})\ \label{arg}
\end{equation}
for $0 \leq i \leq s-1$, where
$$q_i = \frac {\gcd( k-2m_0, k-2m_1,  \dots, k-2m_i)}{
\gcd ( k-2m_0, k-2m_1,  \dots, k-2m_{i+1})}.$$

The model hypersurface $M_P$
associated to  $M$ at $p$ is defined using the normalized leading homogeneous
polynomial,

\begin{equation}\label{2.11} 
M_P = \left\{(z,w) \in \CC{2}:\,\,
 v  = \sum_{j=1}^{k-1}
a_j z^j\bar z^{k-j} \right\}. 
\end{equation}

The following transformations preserve the model:
\begin{equation}
\label{h1}
 z^* =   \delta e^{i\theta} z,\ \  \ \ \   w^* = \delta^k w,\end{equation}
                           where
$e^{i\theta}$ is a $d$-th root  of unity and $\delta > 0$ for $k$
even or $\delta \in \mathbb R \setminus \{ 0 \}$ for $k$ odd.

Next we can apply an analogous strategy as described for the Chern-Moser normal form and subsequent constructions by Beloshapka in Section 2. In particular, the crucial tool is again the Chern-Moser operator, which captures the leading linear part of the transformation formula.
Now it takes the form
$$L(f,g) =  Re \{ig(z,u+i\pz) +
2 P_z f(z,u + i\pz)\}.$$ 

It turns out there are three qualitatively different cases:  generic and two exceptional, depending on the form of the kernel of the Chern-Moser operator.

In particular, when the leading polynomial has a circular symmetry, we write
\begin{equation}S_k = \left\lbrace (z,w) \in \mathbb C^2 \ \vert \   v = \vert z \vert^k
\right\rbrace. \end{equation}

It is preserved by the group of transformations of the form
\begin{equation}\label{h2}
z^* =  \frac{ \delta e^{i\theta} z}{(1 + \mu w)^{2/k}}, \quad  
  w^* =\frac{ \delta^k w}{1 + \mu w},
 \end{equation}

                           where $\delta > 0,$ and $\theta, \mu  \in \mathbb R$.
 Its real dimension is equal to
three. 

Another exceptional model is the tubular hypersurface
\begin{equation}T_k = \left\lbrace (z,w) \in \mathbb C^2 \ \vert \
v = \frac1k \left[ (z+\bar z)^k - 2 Re\; z^k \right]  \right\rbrace. \end{equation}

Analysing the kernel and image of the Chern-Moser operator, we arrive at the following results. The first normal form applies to generic models. 

\begin{theorem}
If $e <\frac{k}2$ and $M_P$ is different from $T_k$, then  there exists a formal change of coordinates
%of the form (\ref{formal}), 
such that
in the new coordinates the defining equation satisfies
the following  normal form conditions
%   $F$ is in  normal form if
$$\begin{array} {rl}
N_{j0} & = 0, \ \ \ \ \ j=1,2,\dots,  \\
N_{k-l+j,l} & = 0, \ \ \ \ \ j= 0,1,\dots, \\
N_{2k-2l, 2l} & = 0, \\
(N_{k-1}, P_z) & = 0,  \end{array}$$
where
$$(N_{k-1}, P_z) = \sum_{j=1}^{k-2}N_{j,k-1-j}  (j+1)\bar a_{j+1}.$$
 Such coordinates are determined uniquely, up to a natural action
of the symmetry group of the model, given by scalings
\begin{equation}
\label{dilations}\Lambda(z,w)=(\lambda z,\lambda^k
w),\,\lambda\in\RR{}\setminus\{0\}.
\end{equation}  
\end{theorem}

In the circular case we obtain the following result.

\begin{theorem}
If $e = \frac{k}2$, there exists a  formal change of coordinates such that
in the new coordinates the defining equation satisfies
% the following  normal form conditions
\begin{equation}
\begin{array}{rl} N_{0j} &= 0, \ \ \ \ \ j=0,1,\dots,  \\
N_{e,e+j} &= 0, \ \ \ \ \ j= 0,1,\dots, \\
N_{2e, 2e} & = 0,\\
N_{3e, 3e} & = 0,   \\
N_{2e, 2e-1}& = 0. \label{fo1} \end{array}
\end{equation}
Such normal coordinates
% i.e.those in which the conditions (\ref{fo1}) hold, 
are determined
uniquely up to the action of the local automorphism group (\ref{h2}).
\end{theorem}

When $M_P = T_k$, we obtain a completely analogous result as for generic models,
 with the following normal form conditions:
\begin{equation}
\begin{aligned} N_{0j} &= 0, \ \ \ \ \ j=1,2,\dots,  \\
N_{1,k-1+j} &= 0, \ \ \ \ \ j= 0,1,\dots, \end{aligned}
\label{f2a}
\end{equation}
and
\begin{equation} N_{2k-2, 2} = \re N_{k-2,1} =  \re N_{k, k-1}  = 0.
 \label{f2}
\end{equation}

Again, normal coordinates are determined uniquely up to the action
of the group \eqref{dilations}.

\subsection{Convergence results in the $\CC{2}$ case}

As a well known fact \cite{ber1}, formal maps between finite type hypersurfaces are convergent, so that the normal form in \cite{kolar} solves the holomorphic equivalence problem, even though the normal form itself can be divergent (as was shown by the first author in \cite{kolardiverg}).  However, as was explained in the introduction, a convergent normal form is still fundamentally important as it yields, for example, an explicit description of the moduli space of CR-manifolds under consideration (see, e.g., \cite{generic}). In addition, a convergent normal form is important for understanding the geometry of weakly pseudoconvex domains, which is still far from being understood completely in Several Complex Variables.  

The convergence problem for the formal normal form in \cite{kolar} was addressed in the work \cite{generic}, \cite{cmhyper}  of the second and the third author. The difficulties in studying the convergence can be roughly described as follows. The convergence proof for the normal form by Chern and Moser is
heavily based on the property that  Levi forms at different points
are equivalent. In fact, the geometry and normal forms look
similar at all points. Consequently, normalization conditions for
the normal form at a point $p\in M$ depend analytically on $p$,
and hence, can be translated into systems of certain analytic
ordinary differential equations whose solutions are again
real-analytic, as was explained in Section 2. This is {\em not the case any more}
 for 
Levi degeneracy points of a Levi-nonflat hypersurface in $\CC{2}$, in whose neighborhoods the Levi form does not
have constant rank. Thus geometry at those degenerate points is
fundamentally different from that at Levi-nondegenerate points. 

In the cited work, the authors introduced a general strategy of possible overcoming this difficulty. The case considered is when a certain geometric data, namely, the type of the hypersurface is still constant along at least one smooth real-analytic curve, passing through the reference finite type $k$ point $p$ and transverse to the complex tangent at $p$. Note that existence of such a curve is natural in the spirit of Chern-Moser type convergent normal form, whose normalization conditions do not change along the chain.

% and where the normal 
% which  Indeed,  considering  normal forms of this kind
%\begin{equation}\Label{basic}
%v = \sum_{j,l>0}\Phi_{jl}(u) z^j \bar z^l, \quad (z,w)\in\CC{}\times\CC{}, \quad w=u+iv,
%\end{equation}
% given by conditions on the coefficient functions $\Phi_{jl}(u)$, we find among them  certain nondegeneracy condition imposed on the leading terms $\sum\Phi_{jl}z^j\bar z^l$ with minimal $j+l$ (such as the nonvanishing of $\Phi_{11}$ in the Chern-Moser case). A convergent normal form of the latter kind is only possible when { the type of the hypersurface is finite and constant along the line \eqref{Gamma}}, which explains the natural character of the above  assumption. 
 
An investigation of possible geometric situations in $\CC{2}$ shows that one either has finitely many curves of constant type, or a (unique!) totally real surface passing through $p$, which is transverse to the complex tangent and consists of points of the same type $k$ as $p$. ( In fact, this surface coincides with the whole Levi degeneracy set of $M$.) In each of the two cases, the normalization strategy is different.
 
In the first case, the collection of transverse curves  of constant type $\gamma_1,...,\gamma_s$ becomes a biholomorphic invariant of $(M,p)$ and can be used for the normalization. We call each of the $\gamma_1,...,\gamma_s$ {\em a degenerate chain}.  Straightening one of the degenerate chains and eliminating terms holomorphic in  $z$   yields a hypersurface \eqref{basic} with {\em terms of total degree $<k$ being not present there}. This illuminates the possibility of succeeding in constructing a convergent normal form by a series of explicit transformations, each of which is responsible for normalizing certain geometric data associated with a hypersurface under consideration. After picking and normalizing one of the degenerate  chains, one arrives at a totally different problem: normalizing a triple
$$(M,\gamma,p).$$ 
For that, we develop a formal theory of such triple, and then hope to achieve the formal normalization conditions by convergent transformations. However, surprisingly it turns out that by using the real defining function one still does not get a convergent normal form, while the convergence is obtained by using the {\em complex defining function} \begin{equation}\label{complex}
w=\Theta(z,\bar z,\bar w)
\end{equation}
of a hypersurface (see, e.g., \cite{ber}) instead of the real defining equation 
\begin{equation}\label{real}
v=\Phi(z,\bar z,u).
\end{equation} 
Recall that \eqref{complex} is obtained from \eqref{real} by substituting $u=(w+\bar w)/2,\,v=(w-\bar w)/2i$ and solving for $w$ by using the implicit function theorem (and \eqref{real} is obtained from \eqref{complex} similarly).
Another important tool for convergence is elimination the complex terms $\sum\Phi_{jl}z^j\bar z^l$ with minimal $j+l$. Finally, the following theorem is proved in \cite{cmhyper}:

\begin{theorem}\Label{5}
Let $M\subset\CC{2}$ be a real-analytic hypersurface of
finite type $k\geq 3$ at a point $p\in M$. 
Assume that there exists at least one but finitely many
transverse curves of constant type  in $M$, passing through $p$.
% Let then $\bigl\{\gamma_1,...,\gamma_s\bigr\}$ denotes the finite collection of %degenerate chains in $M$ through $p$. 
Then, for any choice of a transverse curve of constant type (degenerate chain) $\gamma$ through $p$, there exists a biholomorphic map $F:\, (\CC{2},p)\mapsto (\CC{2},0)$,
sending $\gamma$ into the line \eqref{Gamma}
and
 $(M,p)$ into the normal form 
$$ \Bigl\{w=\bar w+2iP(z,\bar z)+\sum_{\alpha,\beta>0,\,\alpha+\beta\geq k}
\Theta_{\alpha\beta}(\bar w)z^\alpha\bar z^\beta\Bigr\},$$
 such that one of the two following cases {\bf (i)} or {\bf (ii)} holds:

\smallskip

\noindent {\bf (i) (circular case)}  The type $k=:2\nu$ is even,  
$$P(z,\bar z)=|z|^k,$$
and the functions $\Theta_{\alpha\beta}$ satisfy
\begin{equation}\label{strongcirc}
%\Phi_{\nu\nu}=0, \quad
%\Phi_{\alpha 0}=0,\,\alpha\geq 0,\quad 
\Theta_{\nu\alpha}=0,\,\alpha \ge \nu,\quad
%\Phi_{\alpha\beta}=0,\,\alpha+\beta\leq k-1,\quad
\im\Theta_{2\nu,2\nu}=\im\Theta_{3\nu,3\nu}=0. 
\end{equation}
The normalizing transformation $F$  is uniquely determined by the restriction of its differential $dF_p$ onto the complex tangent $T^{\CC{}}_p M$ and the restriction of its Hessian $D^2 F_p$ onto the tangent space $T_p\gamma$. In turn, $F$ is unique up to the action of the subgroup 
\begin{equation}
\label{biggroup}\Lambda(z,w)=\left(\frac{\lambda e^{i\theta} z}{(1+rw)^{1/\nu}},
\frac{\lambda^{2\nu}
w}{1+rw}\right),\,\lambda,\theta,r\in\RR{},\,\lambda\neq 0,
\end{equation}  
of the projective group $\mbox{Aut}\,(\mathbb{CP}^2)$.
%Moreover, two germs $(M,p)$ and $(M^*,p^*)$ are biholomorphically equivalent if and only if for some (and then for any) normal forms $(N,0)$ and $(N^*,0)$ of them there exists a linear-fractional map
%\begin{equation}
%\label{biggroup}\Lambda(z,w)=\left(\frac{\lambda e^{i\theta} z}{(1+rw)^{1/\nu}},
%\frac{\lambda^{2\nu}
%w}{1+rw}\right),\,\lambda,\theta,r\in\RR{},\,\lambda\neq 0,
%\end{equation}   
%transforming $(N,0)$ into $(N^*,0)$. 
Moreover, the degenerate chain $\gamma$ through $p$ is unique in the circular case.
% and s given locally by the equation $\{z=0,v=0\}$ in any normal form of $M$ at $p$.

\medskip

\noindent {\bf (ii) (tubular/generic case)} the polynomial $P(z,\bar z)$ has the form
$$P(z,\bar z) = z^\nu\bar z^{k-\nu}+\bar z^\nu z^{k-\nu}+\sum_{\nu+1\leq j\leq k/2}\bigl(a_jz^j\bar z^{k-j}+\bar a_j\bar z^j z^{k-j}\bigr)$$ for some integer $1\leq \nu < k/2$,  and the functions $\Theta_{\alpha\beta}$ satisfy
\begin{equation}\label{stronggeneric}
\begin{aligned}
%\Phi_{\alpha 0}=0,\,\alpha\geq 0,\quad  
\Theta_{\nu\alpha}=0,\,\alpha\geq k-\nu,\quad
% \Phi_{\alpha\beta}=0,\,\alpha+\beta\leq k-1,\quad 
\re\Theta_{2\nu,2k-2\nu}=0 . 
 \end{aligned}
\end{equation}
The normalizing transformation $F$  is uniquely determined by  the initial choice of a degenerate chain $\gamma$ and by the restriction of the differential $dF_p$ onto the complex tangent $T^{\CC{}}_p M$. Moreover, two germs $(M,p)$ and $(M^*,p^*)$ are biholomorphically
equivalent if and only if, for some choice of degenerate chains $\gamma\subset M,\,\gamma^*\subset M^*$ and some (hence any) respective normal forms
$(N,0)$ and $(N^*,0)$, there exists a linear map \eqref{dilations},
transforming $(N,0)$ into $(N^*,0)$.
\end{theorem}

We shall emphasize an important difference in comparison with the Chern-Moser normal form here. Namely that the normal form in Theorem~\ref{5} is not constructed by a Chern-Moser operator at one point. Motivated by this consideration, in general, we expect a convergent normal form construction to be based on power series conditions along certain invariant analytic subsets rather than at one given point.

In the second case, when there exists a two-dimensional surface consisting of points of constant type, the strategy of the convergence proof is very different. In this case, a formal normal form (which is slightly different in \cite{cmhyper} from the one in \cite{kolar}) provides one with a distinguished  direction field in the Levi-degeneracy set $\Sigma$. Namely,  at each point $q\in\Sigma$ one takes the direction of the vector which is transformed into the distinguished direction \eqref{Gamma} in some (and hence any) normal form at $q$. This gives one a distinguished foliation in $\Sigma$, and each of the leaves is called now {\em a degenerate chain}. For reasons similar to the ones described above, having a degenerate chain straightened to a line, makes it possible to arrive to the formal normal form by means of convergent transformations. This finally yields the following theorem:

\begin{theorem}
Let $M\subset\CC{2}$ be a real-analytic hypersurface of
finite type $k\geq 3$ at a point $p\in M$. Assume that the maximal type locus at $p$ has dimension $2$. 
Then there exists a biholomorphic map $F:\, (\CC{2},p)\mapsto (\CC{2},0)$,
which brings $(M,p)$ into a normal form 
$$\Bigl\{ v=P(z,\bar z)
+\sum_{\alpha + \beta \geq k, \, \alpha,\beta>0}
\Phi_{\alpha\beta}(u)z^\alpha\bar z^\beta\Bigr\},$$
 where the polynomial 
 $P(z,\bar z)=\frac{1}{k}\Bigl[(z+\bar z)^k-z^k-\bar z^k\Bigr]$, and $\Phi_{\alpha\beta}$ 
 satisfy 
\begin{equation}\label{strongtube2}
\begin{aligned}
%\Phi_{\alpha 0}=0,\,\alpha\geq 0,\quad 
\Phi_{\alpha 1}=0,\,\alpha \ge k-1,\quad
 \re\Phi_{k,k-1}=\im\Phi_{2k-2,2}=0 .& 
 \end{aligned}
\end{equation}
The normalizing transformation $F$  is uniquely determined by the restriction of its differential $dF_p$ onto the complex tangent $T^{\CC{}}_p M$. Moreover, two germs $(M,p)$ and $(M^*,p^*)$ are biholomorphically equivalent if and only if, for some  (and then for any)  normal forms $(N,0)$ and $(N^*,0)$ of them, there exists a linear map \eqref{dilations}
transforming $(N,0)$ into $(N^*,0)$. The Levi degeneracy set $\Sigma$ of $M$, which in this case is a smooth real-analytic totally real surface in $\CC{2}$ transverse to the complex tangent $T^{\CC{}}_p M$  is canonically foliated by distinguished biholomorphically invariant curves, called degenerate chains, where the chain through p is locally given by $\bigl\{z = 0,\,v=0\bigr\}$ in any normal form at p.
\end{theorem}

%We emphasize that the approach used in \cite{cmhyper} can be generalized for a wide range of finite type hypersurface where the constancy of certain geometric data is present.

%\begin{remark}
%The outlined above proof of Theorem 6 can be used to simplify the original convergence proof of Chern and Moser, as was mentioned in Section 2. One should use the formal normal form of Chern and Moser to obtain a distinguished direction field in the bundle of $1$-jets of curves transverse to the complex tangent at every point (by choosing for every point in the jet bundle the unique direction transformed into the prolonged direction \eqref{Gamma} by the jet prolongation of some, and then any, transformation bringing $M$ to a normal form). After that, the canonical projection of the leaves of the respective foliation provides us with chains in $M$. After straitening a chain, we repeat the sequence of explicit Chern-Moser transformation and arrive to the desired normal form.
%\end{remark}

\subsection{Higher-dimensional case}

In complex dimensions higher that two, there exist several different generalizations of the concept of finite type. 
In the normal form context, it turns out that the Chern - Moser operator can  
be  extended to hypersurfaces of finite Catlin multitype.
It formalizes  the notions of model hypersurface
and invariantly defined weights, which are both essential for Chern-Moser theory.

We recall that multitype in the sense of Catlin is a fundamental CR invariant which 
Catlin introduced  to prove subelliptic estimates on pseudoconvex domains (see \cite{catlin}). 
 On the other hand, in our context there are  no pseudoconvexity assumptions
(multitype was considered in the general case in \cite{kol10}), similarly as the   Chern-Moser theory allows  model hyperquadrics of all signatures.

Let  now $M \subseteq \mathbb C^{n+1}$ be a smooth hypersurface,
and $p \in M $ be a  point 
on $M$. 
We will consider
local holomorphic coordinates $(z,w)$ centered at $p$,
where $z =(z_1, z_2, ..., z_n)$ and  $z_j = x_j + iy_j$,
$w=u+iv$. The hyperplane $\{ v=0 \}$ is assumed to be tangent to
$M$ at $p$, hence  $M$  is described near $p$ as the graph of a uniquely
determined real valued function
\begin{equation} v = \psi(z_1,\dots, z_n,  \bar z_1,\dots,\bar z_n,  u), \ d\psi(p) \neq 0.
\label{vp1}
\end{equation}

The definition of multitype assigns   rational  weights  to the variables
$w, z_1, \dots z_n$. To begin with, the
variables $w$, $u$ and $v$ are given weight one, reflecting our choice
of  complex normal variables.
The weights of complex tangential variables $(z_1, \dots, z_n)$ are determined as follows. 
By a weight we will understand an $n$-tuple of nonnegative
 rational numbers $\La = (\la_1, ...,
\la_n)$, where $0 \leq\la_j\leq \frac12$, and $\la_j \ge
\la_{j+1}$. Now, we give

\begin{definition} A weight $\La$ will be called distinguished for $M$ if there exist
local holomorphic coordinates $(z,w)$ in which the defining equation of $M$ takes form
\begin{equation} v = P\zz + o_{\La}(1),
\label{1}
\end{equation}
where $P\zz$ is a nonzero $\La$ - homogeneous polynomial of
weighted degree $1$ without pluriharmonic terms, and $o_{\La}(1)$
denotes a smooth function whose derivatives of weighted order less than or equal to
one vanish.
\end{definition}

%The fact that distinguished weights do exist follows from \eqref{fifi}.
%For these coordinates $(z,w),$ we have
%$$\Lambda=(\dfrac{1}{m}, \dots, \dfrac{1}{m}).$$

We will consider the lexicographic order on the set
of weights,  defined in a usual way:
$(\a_{1},\ldots,\a_{n})<(\b_{1},\ldots,\b_{n})$
whenever for some $1\le k\le n$,
$\a_{j}=\b_{j}$ for $j< k$ but $\a_{k}<\b_{k}$.

We now recall the definition of  Catlin multitype.

\begin{definition}
Let  $\Lambda_M = (\mu_1, \dots, \mu_n)$  be the infimum of all possible
distinguished weights  $\Lambda$ with respect to the lexicographic order.
The multitype of $M$ at $p$ is defined to be the $n$-tuple $$(m_1,
m_2, \dots, m_n),$$ where
$$m_j = \begin{cases}   \frac1{\mu_j} \ \  {\text{ if}} \ \  \mu_j \neq 0\\
  \infty \ \ {\text{ if}} \ \   \mu_j = 0.
\end{cases} $$
\end{definition}

Furthermore, if none of the $m_j$ is
infinity, we say that $M$ is of {\it finite multitype at $p$}.

Let us now assume that $ p \in M $ is a point of finite  Catlin multitype
 $(m_1, \dots, m_n)$.
 
 It is not difficult to see (\cite{catlin}) that 
 the infimum in the above definition is in fact attained, and  one can find coordinates $(z_1, \dots, z_n,,w)$ with weight of  $z_j$ equal to  $\mu_j =  \frac1{m_j}$,
 weight of  $w$ equal to 1 such that
 $M$ is described by
\begin{equation}\Label{basic}
\Im w = P(z,\bar z) + F(z,\bar z,\Re w),
\end{equation}
where $P$ is a weighted homogeneous polynomial of weighted degree  $1$,
and $F$ has Taylor
expansion  with terms of weighted degree $>1$.

 Assume that $M$ is given by 
 (\ref{basic}), and the associated model hypersurface
\begin{equation}\Label{model2}
M_{P}:=\{\Im w = P(z,\bar z)\}
\end{equation}
%\eqref{model}
is holomorphically nondegenerate.
%The following theorem shows that under these assumptions,  we can describe the  %infinitesimal automorphisms of the model hypersurface.

 The Chern-Moser operator is now defined as
 \begin{equation}\label{chern}
L(f,g)= \Re \left\lbrace ig(z,u+iP(z,\bar z)) +
2\sum_{j=1}^{n} P_{z_j}(z, \bar z)  f^j(z,u + i\pz)\right\rbrace.
\end{equation}

The following two results were proved in \cite{KMZ}.  Here by   $E$ we denote the set 
% all rational numbers in the interval $(0,1)$ that can be written
%as 
\begin{equation}
E = \left\lbrace  \sum_{j=1}^n k_j \mu_j;\ \ \ k_j \in \mathbb N \cup \{-1 \} \right\rbrace \cap (0,1).
\end{equation}
%combinations of the weights $\mu_1, \dots, \mu_n$ with coefficients in $\{-1\} %\cup \mathbb N$. 

\begin{theorem}
 Let $M_{P}$ be given by \eqref{model2}.
The   Lie algebra of infinitesimal automorphisms
$\5g=\aut (M_{P},0)$ of $M_{P}$ admits the  weighted grading given
by
\begin{equation}
\5g = \5g_{-1} \oplus \bigoplus_{j=1}^{n}\5g_{- \mu_j} \oplus \5g_{0}
\oplus \bigoplus_{\eta \in E}\5g_{\eta}
%\oplus \bigoplus_{j=1}^{n}\5g'_{1- \mu_j}
 \oplus \5g_{1}.
 \label{muth}
\end{equation}
\end{theorem}

We will denote by $(f_1, f_2, \dots, f_n, g)$ the components 
of an automorphism of $M$, as in Section 2,  and by $\vert \alpha \vert_{\Lambda_M}$
the weighted length of a multiindex $\alpha$.

\begin{theorem}

The automorphisms of $M$ at $p$  are uniquely determined jointly by:
%the following partial derivatives
\begin{enumerate}

\item the complex tangential derivatives
$\dfrac{\partial^{\vert \alpha \vert}f_j}{ \partial z^{\alpha}}$ for
$\vert \alpha \vert_{\Lambda_M} \leq 1-\mu_n$;

\item the first and second order normal derivatives
$\dfrac{\partial f_j}{\partial
w}$ for $j= 1, \dots, n$,  $\dfrac{\partial g}{\partial
w},$ $\dfrac{\partial^2 g}{\partial
w^2}.$

\end{enumerate}
 \end{theorem} 
 
 One of the cases of finite Catlin multitype was considered in the work \cite{generic}. There the authors considered the case of  generic Levi degeneracy points, as introduced by Webster \cite{webstergeneric}, i.e. points $p\in M$ where the Levi determinant vanishes but its differential
does not and the set of Levi degenerate points of $M$ (which is then a
smooth codimension-one submanifold of M at $p$) is transverse to the
Levi null space (which is then one dimensional) at $p$. Note that a point of generic Levi
degeneracy is ``stable'' in the sense that it cannot be removed by
a small perturbation. The multitype here is $\left(\frac{1}{2},...,\frac{1}{2},\frac{1}{3}\right).$ A different {\em formal normal form} for
generic Levi degeneracy points, whose convergence remains unknown,
was constructed by Ebenfelt \cite{ebenfeltjdg} (where the author did not use the strategy of homological operators). In \cite{generic}, the authors use the strategy of Chern-Moser operators to obtain a formal normal form with the convergence property. The proof of the convergence there is somehow similar to that for Theorem 6, that is why we omit discussing the formulation and the proof here  (in fact, the paper \cite{generic} preceded the paper \cite{cmhyper} where Theorem 6 was proved, however, we chose Theorem 6 for the detailed description as it is technically simpler).  

\subsection{Normal forms in the infinite type setting}

Local geometry of infinite type hypersurfaces turns out to be more complicated and in many ways strikingly different from the finite type situation.
So far the normal form approach has been applied only in some special  cases
(\cite{elz}, \cite{kl}).

Here we consider the class of 
ruled infinite type hypersurfaces in $\mathbb C^2$, which  we    
denote  by $\mathcal{A}$. Such hypersurfaces  play the role of model hypersurfaces for 
the more general class of $1$-nonminimal hypersurfaces. 

% Hypersurfaces in  $\mathcal{A}$  can be given by an equation of the form 
% \begin{equation}
% \Im w = (\Re w ) A(z,\bar z). 
% \label{aa}
% \end{equation}
% and with a slight abuse of notation, we identify the germ of the hypersurface 
% with the germ of $A$.

Let us consider a real analytic ruled hypersurface
 $M \subseteq \mathbb{C}^2$. We will assume that $p\in M$ is a 
point of infinite type on $M$. It implies that there is a complex hypersurface
$E$ passing through $p$ and contained in $M$. The real lines which form $M$ then intersect $E$ transversally. 

 We  choose local coordinates $(z,w)$ in $\mathbb C^2$ in which $p=0$ and such that locally around $p$, 
$E= \left\{ w=0\right\}$.
In such coordinates, $M$ is described by 
\begin{equation}
        \label{e:ruled1} \Im w = (\Re w) \, A(z,\bar z).
\end{equation}
 With a slight abuse of notation, we identify the germ of the hypersurface 
with the germ of $A$.
We would  
like to decide whether two
such hypersurfaces are biholomorphically equivalent.

Let us again order points in $\mathbb N^2$  lexicographically. We write $(a,b) << (c,d)$ if
either $b< d$ or if $b=d$ and $a\leq c$. For 
such an $A(z,\chi) = \sum_{\alpha,\beta} A_{(\alpha,\beta)} z^\alpha \chi^\beta \in\mathcal{A}$, we define
\begin{equation}
        \label{e:initialterm} \initial A = \min \{(\alpha,\beta)\in \N^2 \colon A_{(\alpha,\beta)} \neq 0 \}, 
        \text{ and  }
        \initialterm A = A_{\initial A}
\end{equation}
where the minimum is with respect to the lexicographic order.
It is not difficult to see that the preceding definition is indeed independent of 
the choice of coordinates, and thus gives rise to a local invariant of $(M,0)$. The normalization of 
$A\in\mathcal{A}$ is now done by a change of coordinates of the form $(z',w')=(\psi(z),w)$ and reduces
the equivalence problem in $\mathcal{A}$ to a linear question:
\begin{theorem}\label{lem:normalization}
        Let $A\in \mathcal{A}$ be given. 
        Then there exists a holomorphic change of coordinates preserving the form \eqref{e:ruled1} such that in the new coordinates 
        $(z',w') = (\psi(z), w)$, the hypersurface corresponding to $A$ is given by $A' (z',\chi')$ 
       and satisfies
        \begin{equation}
                \label{e:normalizationconditions} 
                A'(z',0) = A'(0,\chi) = 0, \quad A'_{\initial A + (0,k)} = 0, \quad k>0, \quad \initialterm A' = \begin{cases} 1 & \initial A \neq (n,n), \\
\pm 1  & \initial A = (n,n) . \end{cases}
        \end{equation}
        We shall say that $A'$ is in normal form if it satisfies \eqref{e:normalizationconditions}. If $A$ 
        and $A'$ are both in normal form, and $H=(f(z,w), g(z,w))$ is a germ of a biholomorphism which transforms
        $A$ into $A'$, then $H=L\circ T$, where $L$ is a linear map of the form $(z,w)\mapsto (\lambda z, w)$ with
         $|\lambda|=1$ and
        $T$ is an automorphism of the hypersurface corresponding to $A$.
\end{theorem}

The next issue addresses in \cite{kl} using a normal form approach is the description of all ruled hypersurfaces which admit automorphisms nonlinear in normal coordinates.  
The results  show that there exist essentially 3 families of ruled hypersurfaces which possess 
nonlinear automorphisms,
which we denote  by $\mathcal{A}_1$, $\mathcal{A}_2$, and $\mathcal{A}_3$: 

\begin{enumerate}
        \item $\mathcal{A}_1$ consists of hypersurfaces which are preimages of 
         Levi-nondegenerate hyperquadrics in $\C^{2k-1}$, given 
        by $\im \eta = \re \left(\sum_{j=1}^{k-1} p_j \zeta_j \overline{\xi_j} \right)$ (where
        $p_j= \overline{p_{k-j}}$), under
        the map
         \[\eta = w, \quad \zeta_j = z^j w^{\left(1-\frac{j}{k}\right)}, \quad 
                \xi_j = z^{k-j}w^{\frac{j}{k}}.\]
                The elements of  $\mathcal{A}_1$  are {\em rational blowups of hyperquadrics}. A hypersurface in $\mathcal{A}_1$ is completely determined by $k$ and $p_1,\dots,p_{(k-1)/2}$. 
                
        \item $\mathcal{A}_2$ consists of hypersurfaces which are preimages of a circular 
        finite type hypersurface in $\C^2$ given by $\im \eta = |\zeta|^{2\ell}$ under a
        map of the form \[ \zeta = z w^{\frac{1}{2\ell} + iT}, \quad \eta = \pm w,\]
        where $T\in\R$.
      Elements of  $\mathcal{A}_2$ are {\em transcendental blowups of the ball} if $T\neq 0$. 
      
        \item $\mathcal{A}_3$ can be  described as consisting of {\em tube-like
        hypersurfaces}.
        %Hypersurfaces here are given by an equation $\im w = \re w A(z,\bar z)$. 
        Hypersurfaces in 
        $\mathcal{A}_3$ are uniquely identified by the property that they possess an infinitesimal 
				automorphism of the form $ w^2 \dz $. Every such hypersurface is uniquely determined 
          by a  part of the Taylor expansion of its defining function $A(z,\bar z)$, namely, 
        $A(z,0) $. In contrast to 
 				the $\mathcal{A}_1$ and $\mathcal{A}_2$, the parameter space here is infinite dimensional. 
\end{enumerate} 

In order to formulate the main result, let 
us  say that a hypersurface $\im \tilde w = \re \tilde w \tilde A (z,\bar z)$ is a root
of a hypersurface $\im w = \re w A(z,\bar z) $ if it is obtained from the latter by 
a modification of the form $\tilde w^k = w$. 

\begin{theorem}\label{thm:mainnonlinear} Let $M$ be a ruled hypersurface in $\mathbb C^2$, 
        $p\in M$ of infinite type. If $(M,p)$  allows an automorphism which 
        is nonlinear in normal coordinates, then $M$ has a root 
        in $\mathcal{A}_1$, $\mathcal{A}_2$, or $\mathcal{A}_3$.
\end{theorem}

\begin{remark} It is important to emphasize that, in the $1$-nonminimal case, a formal normal form would solve the equivalence problem, as follows from the convergence result \cite{jl2} of Juhlin and Lamel. However, in the general infinite type case, one can not hope for formal normal forms to solve the holomorphic equivalence problem, since formally equivalent real hypersurface can be still inequivalent holomorphically, as was proved in the work \cite{divergence} of Shafikov and the second author.
\end{remark}

\section{Symmetry preserving normal forms}

The  normal form constructions of
\cite{chern}, \cite{kolar} have two rather unpleasant features, 
which may substantially reduce their applicability.
First, the normal forms
are uniquely determined only up to the action of the symmetry group of the model.
In order to verify the
local biholomorphic equivalence it is in general not sufficient to compose
the normal form with an element of the symmetry group. 
This may actually produce power series which is not in normal form, and one
has to renormalize.
In practice, this may be impossible to carry out.

The second problem comes from the fact that the above normal forms do not in general
respect symmetries of the hypersurface.
 For analysis on a  hypersurface which is known to admit a symmetry,
it would  be  desirable to use  a normalization  which reflects
this symmetry.

 This problem was first addressed  by N. K.
Stanton, who considered so called
 rigid hypersurfaces  in $\mathbb C^2$,
 and constructed a normal form which respects the underlying translational symmetry.
 In fact, the results of \cite{stanton} describe
  all transformations
 preserving this normal form, and give a complete
 classification of rigid hypersurfaces, provided that the model does not admit
 an additional symmetry.

The approach of Stanton was generalized in \cite{eks}. It is based on  the fact that
automorphism groups are essential  geometric invariants,
which should be taken into account before any attempt to compute
higher order invariants.

Let us first consider a hypersurface $M$, which has a one dimensional symmetry group, consisting of transversal
symmetries. A complete normal form construction for such hypersurfaces follows
from  the work of Stanton.

We will use  local holomorphic coordinates in which the corresponding vector field is of the form
$$X = \frac{\partial }{\partial w}.$$
The defining equation of $M$ is then  independent of $u$, hence
$$v = F(z, \bar z).$$
 Since we assume that $M$ admits no other automorphism, it follows that
all transformations which
preserve such a form have to preserve the corresponding vector
field, up to a real multiple. Indeed, if $X$ maps to a different field, then
this field gives an additional symmetry of $M$ 
(see also Proposition 3.1. in \cite{stanton}).
We can therefore immediately  verify that such transformations can be written as

\begin{equation}\begin{aligned}
  z^* &=  f(z)\\
  w^* &= cw + g(z), \\
\label{wz}
\end{aligned} \end{equation}
where $c\in \mathbb R^*$ and  $g(0) = f(0) = g'(0) = 0$, $f'(0) \neq 0$.

In order to state Stanton's normal form result, we will 
 use the expansion of $F$ in terms of $z, \bar z$,
$$F(z, \bar z) = P_k(z, \bar z) + \sum_{j,l}A_{jl}z^j \bar z^l,$$
where $P_k$ if of the form \eqref{pzz}, and satisfies the normalization condition
\begin{equation}
 a_e =1,
\label{al1}
\end{equation}
where e is as defined in Section 3.
The normal form conditions are
\begin{equation}
 A_{0,l}=A_{e,l}=0,
\label{snf}
\end{equation}
for all $l = 1,2,\dots$

Stanton proved the following result (Theorem 1.7. and Proposition 3.1 in \cite{stanton}).

\begin{theorem}
There exists a transformation of the form (\ref{wz}) which takes $M$ into normal form.
If $P_k$ is rotational, the normal form
is determined uniquely
up to a rotation in $z$ and a weighted dilation, i.e. a transformation of the form
\begin{equation}
 z^* = \la e^{i\theta} z, \ \ \ \ \ w^* = \pm \la^k w.
\end{equation}
If $P_k$ is not rotational, the normal form
is determined uniquely
up to  a weighted dilation.
\end{theorem}
%Stanton's normal form is unique up to weighted dilations.
Since the action of rotations and dilations on the defining equation is completely
straightforward, the remaining real parameter(s) can be
easily fixed (see \cite{eks}).

We now consider  a hypersurface $M$ which has  a one dimensional
algebra of nontransversal symmetries.

In the following, we will use a standard weight assignment,
the complex tangential variables $z,x,y $ are given weight
 one and  $w$ and $u$ weight $k$.

Consider local holomorphic coordinates in which the infinitesimal CR automorphism
takes the form
% By straightening the corresponding holomorphic vector
%field, it is taken into the form
$$Y= i \frac{\partial }{\partial z}.$$
In such coordinates, the defining equation is

\begin{equation} v =  G(x,u), \label{vg}\end{equation}
where
\begin{equation}
G(x,u)=  x^k + \sum_{j,l=1}^{\infty} A_{jl} x^j u^l, \label{ff}
\end{equation}
and  the sum on the right contains terms of weight bigger than $k$.
We use also the partial expansion of $G$,
$$G(x,u) = x^k + \sum_{j=0}^{\infty}
X_j(u) x^j,$$
where
\begin{equation}
X_j(u) = \sum_{l}A_{jl} u^l.
\end{equation}

By the same reasoning as in the previous case, since $M$ admits no other automorphisms, the only transformations
preserving form (\ref{ff}) have to preserve the vector field, up to a real multiple
(otherwise, the field into which
$Y$ maps gives an additional symmetry of $M$).
Thus the transformations are of the form
\begin{equation}
z^* = cz + \psi(w), \ \ \ \ \  w^* = \phi(w), \label{ps}
\end{equation}
\\[2mm]
where
$c\in \mathbb R^*$ and  $\psi(0) = \phi(0) = 0$, $\phi'(0) \in \mathbb R^*$.
\begin{definition}{
$M$ is in normal form if the defining equation  has form
(\ref{vg}), and $G$ satisfies
\begin{equation}
 X_{0} =X_{k-1}=
X_{k}=  X_{2k-1} = 0. \label{a}\end{equation} }
\end{definition}
In particular, when $k=2$ and $M$ is Levi nondegenerate, the
normal form conditions are
\begin{equation}
 X_{0} =X_{1}=
X_{2}=  X_{3} = 0. \label{a3}\end{equation}
\\[2mm]

\begin{theorem}
{There exists a unique formal transformation of the form
\begin{equation}
z^* = z + \sum_{j=1}^{\infty} c_j w^j, \ \ \ \ \  w^* = w +
\sum_{j=2}^{\infty} d_j w^j \label{ps2}
\end{equation}
which takes $M$ into normal form.}
\end{theorem}

As it was shown in \cite{eks} in the case of a 2-dimensional
Abelian symmetry algebra $\mathfrak g$ of $M$ there exist holomorphic
coordinates $(z,w)$ where the two generators of $\mathfrak g$ take
the form
$$X=\frac{\partial}{\partial z},\quad Y=\frac{\partial}{\partial w}.$$

and the equation of the germ $(M,0)$ is
% \sideremark{perhaps we need not separate the 2k-1 term here.}
$$v=G(y)=ay^k+by^m+\sum \limits_{j=m+1}^\infty G_j y^j ,$$
with $2\le k< m\neq 2k-1$ and $a,b \in \R \setminus \{0\}$.
The only possible coordinate change preserving the symmetry algebra $\mathfrak g$ of
$M$ and the complex tangent $T_0^CM=\{w=0\}$ is linear: $ z \mapsto \alpha z +\beta
w,\quad w \mapsto \gamma w, $ where $\alpha,\, \beta, \, \gamma \in \mathbb R.$
Appropriately and uniquely choosing $\alpha, \beta$ and $\gamma$ we obtain the tubular normal form of $M$: 
$$v=G(y)=y^k + \varepsilon y^m+ \sum_{j=m+1}^\infty G_j y^j ,$$
with $2\le k< m\neq 2k-1$ and $G_{2k-1}=0$. Here 
$$\varepsilon=\begin{cases}1 \text{ if } $m$ \text{ is odd and}\\\pm 1 \text{ if } $m$ \text{ is even.}\end{cases}.$$

Notice that the integer $m$, together with the type $k$, is a
biholomorphic invariant of $(M,0)$. The hypersurfaces
$$v=y^k + \varepsilon y^m$$
can be considered as the models for this type of symmetry.

\section{Open problems}

In this section we would like to give a list of open problems which we find particularly interesting for the development of the field and for applications.

The first problem concerns Beloshapka's strategy for normal forms of Levi-nondegenerate hypersurfaces. We recall that in any CR-dimension $n>0$ and CR-codimension $k>0$ an important object is the moduli space $\mathcal M(n,k)$ of the respective models which for some values of $n,k$ has positive dimension. The latter can be also interpreted as the fact that Levi-Tanaka algebras at different points in $M$ are non-isomorphic. This is an effect which, for example, is ultimately excluded from the consideration in the celebrated Tanaka's theory of prolongation in Cartan geometry or in the work of the school of Cap  (see Tanaka \cite{tanaka}, Cap and Schichl \cite{cap-schichl}).   In all known instances, the authors excluded this case from consideration in normal form constructions as well, since it looks difficult to handle the non-uniformity of the Chern-Moser operator when changing the base point. That is, the equivalence problem in this case is completely open, even though, for example, the class of  Levi-nondegenerate codimension $2$ submanifolds in $\CC{5}$ is already within this case.   In the light of this, we formulate

\begin{problem} Construct a complete normal form for Levi-nondegenerate codimension $2$ submanifolds in $\CC{N},\,N\geq 3$. 
\end{problem}

The second problem is related to infinite type hypersurfaces. Besides the normal form in \cite{kl}, no constructions are known in this setting. The main difficulty here is the absence of polynomial models for this class of hypersurfaces (see, e.g., the discussion in the Introduction in \cite{nonminimalODE}). An additional difficulty here follows from the remark in the end of Section 3. 

\begin{problem} Construct a normal form for infinite type hypersurfaces in $\CC{2}$, and investigate the convergence issue.
\end{problem}

This problem remains one of the most intriguing in CR-geometry. A possible approach here is given by the CR -- DS technique (see \cite{nonminimalODE}) where an infinite type hypersurface is replaced by a second order singular ODE, which already {\em do admit} certain polynomial models (known as {\em Poincare-Dulac normal forms} or {\em Birkhoff normal forms}, see \cite{ilyashenko}). 
We expect such a normal form would involve certain types of resonances, and the convergence
issue will depend on the presence of some analogue of small divisors.

The third problem occurs in the important pseudoconvex case:

\begin{problem} Construct a convergent normal form for an isolated weakly pseudoconvex point in a hypersurface $M\subset\CC{2}$.
\end{problem}

The main difficulty here is that the known convergent constructions \cite{generic},\cite{cmhyper} are based on the presence of an entire curve passing through the reference finite type point and having certain geometric data constant along it. No such curve exists at an isolated weakly pseudoconvex point. In relation to this, we state the following

\begin{problem} Is it true that at any isolated  weakly pseudoconvex point $p$ in a real-analytic hypersurface $M\subset\CC{2}$ there exists  a smooth real-analytic curve $\gamma\subset M,\,\gamma\ni p$ which is a Chern-Moser chain at any point $q\in\gamma$ with $q\neq p$? Can that curve be chosen to be smooth?
\end{problem}

In the case of presence of such a smooth distinguished curve, it can be straitened and the rest
of the normalization procedure would yield the normalization of a triple $ (M;\gamma 
; p)$, as described in
Section 3.

\end{document}